\documentclass[12pt]{amsart}
\usepackage{amscd,amsmath}
\usepackage{latexsym}
\usepackage{amssymb}
\usepackage{amscd}
\usepackage{amsmath}

\newtheorem{lemma}{Lemma}[section]
\newtheorem{theorem}{Theorem}[section]

\newtheorem{proposition}{Proposition}[section]

\newtheorem{remark}{Remark}[section]

\newcommand{\Z}{\mathbb{Z}}
\newcommand{\CP}{\mathbb{CP}}
\newcommand{\C}{\mathbb{C}}
\newcommand{\R}{\mathbb{R}}
\newcommand{\RP}{\mathbb{RP}}

\newcommand{\md}[1]{\ensuremath{\overline{M}_{k}(#1,\beta)}}

\newcommand{\dm}[1]{\ensuremath{\overline{M}_{#1}}}

\begin{document}
\title[Real Aspects of the Moduli Space ]%
{Real aspects of the
 moduli space of genus zero stable maps}

\author{Seongchun Kwon }


\maketitle

\vspace{1mm}

\begin{abstract}
We show that the moduli space of genus zero stable maps is a real
projective variety if the target space is a smooth convex real
projective variety. We show that evaluation maps, forgetful maps are
real morphisms. We analyze the real part of the moduli space.
\end{abstract}

\section{Introduction}\label{s;intro}
\vspace{1mm}

We call a projective variety $V$ as a \emph{real} projective variety
if $V$ has an anti-holomorphic involution $\tau$ on the set of
complex points $V(\C)$. By a \emph{real structure} on $V$, we mean
an anti-holomorphic involution $\tau$. The \emph{real part} of $(V,
\tau)$ is the locus which is fixed by $\tau$.

\par In the following paragraph, readers can find the definitions of the moduli space of stable maps and
various maps defined on it in \cite{cak}.

\par Let's assume that $X$ is a convex real projective variety. We show the following:

\begin{itemize}
\item The moduli space $\md{X}$ of $k$-pointed genus 0 stable maps is a
real projective variety.
\item Let $\dm{k}$ be the Deligne-Mumford moduli space of $k$-pointed
genus 0 curves. The forgetful maps $F_{n}: \overline{M}_{n}(X,
\beta) \rightarrow \overline{M}_{n-1}(X, \beta)$, $F: \md{X}
\rightarrow \dm{k}$ are real morphisms(i.e., morphisms which commute
with the anti-holomorphic involution on the domain and that on the
target), where $n \geq 1$, $k \geq 3$.
\item Let $ev_{i}: \md{X} \rightarrow X$ be the $i$-th evaluation
map. Then, $ev_{i}$ is a real map.
\item Let $\CP^n$ have the real structure from the complex conjugation
involution. Let $X$ be a real projective variety such that the
imbedding $i$ which decides the real structure on $X$ has a
non-empty intersection with $\RP^n \subset \CP^n$. Let $M_{k}(X,
\beta)$, $k \geq 3$, be the moduli space of $k$-pointed genus 0
stable maps with a smooth domain curve. Then, each point in the real
part $M_{k}(X, \beta)^{re}$ of the moduli space represents a real
stable map having marked points on the real part of the domain
curve.
\end{itemize}

This paper is organized as follows. In Sec.\ref{s;aspect}, we show
that the moduli space is a real projective variety. We prove that
straightforwardly, based on the explicit tangent space splitting
calculation. In Sec.\ref{s;property}, we show that the forgetful
maps, the evaluation map are real maps. Also, we do real part
analysis when $k \geq 3$.
\par Theorem \ref{t;maint} is the main theorem in this paper. Real part analysis
done in sec.\ref{s;property} shows that the studies of the
intersection theoretic properties on the real part of $\md{X}$ are
important for real enumerative applications. The main results of
this paper are similar to those in \cite{kwont}. The main Theorem in
\cite{kwont} and in this paper is based on the proof which shows
that the defined involution on the \emph{complex} moduli space of
stable maps is an anti-holomorphic involution, when the target space
$X$ is a real convex projective variety. However, the practical
methods of proofs are different.

\par The real version of the Gromov-Witten invariants are defined in
\cite{kwon}, when the target space is a rational projective surface.
Different from the Gromov-Witten invariant defined on $\md{X}$, the
real version of the Gromov-Witten invariants are local invariants.

\par Relevant theory which considers the global minimum bound of real enumerative problems
has been developed by J-Y Welschinger in \cite{wel4}, \cite{wel3}.
The Gromov-Witten invariant in the real world with Quantum Schubert
calculus has been widely studied by F. Sottile. See
~\cite{sot1}, ~\cite{sot2}, \\
~\cite{sot3}, \cite{sot4},
~\cite{sot5}.

\newpage

\par \textbf{Convention:}
\begin{itemize}
\item The real structure on $\CP^1$, defined by
the standard complex conjugation involution, the real structure of
the target space $X$ will be always denoted by $s$, $t$
respectively.
\item Let $C$ be an arithmetic genus 0 curve. Let $\pi: \widetilde{C}:= \CP^{1}_{1} \cup \ldots \cup \CP^{1}_{l}
\rightarrow C$ be a normalization map. We will denote the
irreducible component in $\widetilde{C}$ by $\CP^{1}_{q(p)}$, either
if it contains
 $\pi^{-1}(p)$ where $p$ is a non-singular point in $C$, or if it contains $p$ where $p$ is any point in
 $\widetilde{C}$.
\end{itemize}

\section{Real Aspects of the moduli spaces}\label{s;aspect}

The following Lemma is well-known. See
~\cite[2.3]{gam},~\cite[sec.10]{foh}, \cite[4.1]{liu}. However, the
author couldn't find the proof. Thus, we include the proof.
\par The tangent space calculations done in Lemma \ref{L:D-M}, Theorem
\ref{p;tansp}, are from repeated $K$-group calculations of vector
spaces based on the simple homological algebraic fact( cf.
Proposition 2.11. in \cite{aam}). The $K$-group we consider is the
Grothendieck group of $K_{0}(\text{point})$ because the tangent
space is calculated pointwise. The way to express the tangent space
of the Deligne-Mumford moduli space $\overline{M}_{k}$ in the proof
of Lemma \ref{L:D-M} is somewhat different from the conventional
one. However, they are equivalent in the $K$-theoretic point of
view. The alternative expression is taken because it makes us to
easily relate the underlying real structure of the pointed curve
with the anti-holomorphic structure on the Deligne-Mumford moduli
space. We include the details of the proof for this alternative
expression. Note that each element in $\overline{M}_{k} \setminus
M_{k}$ represents a singular curve having only nodal singularities.
Different from the Deligne-Mumford moduli space of pointed higher
genus curves, singular curves in the Deligne-Mumford moduli space of
pointed genus zero curves are trees. Therefore, the number of
singular points(i.e., gluing points) is exactly one less than the
number of irreducible components.

\begin{lemma}\label{L:D-M}
 Let $\mathbf{c}:=[(C,a_{1},\ldots, a_{k})]$ be a point in
$\dm{k}$. Let $\pi: \widetilde{C}:= \CP^{1}_{1} \cup \ldots \cup
\CP^{1}_{l}\rightarrow C$ be a normalization map, where
$\CP^{1}_{i}$ is biholomorphic to $\CP^1$. Let
 $g_{1}, \ldots, g_{r}$, $r=l-1$, be singular points in $C$. Let's denote two
points in $\widetilde{C}$ corresponding to $\pi^{-1}(g_{i})$ by
$g_{i}^{1}, g_{i}^{2}$.
\par Let $\overline{\mathbf{c}}$ be the point in $\dm{k}$
represented by the pointed curve $(\overline{C}, s(a_{1}), \ldots,
s(a_{k}))$ which satisfies the following:\\
Let $\overline{\pi}: \widetilde{\overline{C}}:= \CP^{1}_{1} \cup
\ldots \cup \CP^{1}_{l}\rightarrow \overline{C}$ be a normalization
map.
\begin{itemize}
\item If $\widetilde{g}_{1}, \ldots, \widetilde{g}_{r}$ are singular
points on $\overline{C}$, then, $\widetilde{g}_{i}^{1} \in
\CP^{1}_{q(g_{i}^{1})}$, $\widetilde{g}_{i}^{2} \in
\CP^{1}_{q(g_{i}^{2})}$ in $\overline{\pi}^{-1}(\widetilde{g}_{i})$
are $s(g_{i}^{1}), s(g_{i}^{2})$.
\item $s(a_{i})$ is the point in $\CP^{1}_{q(a_{i})}$ conjugate to $a_{i}$
by the real structure $s$ on $\CP^1_{q(a_{i})}$.
\end{itemize}

The involution, $I: \mathbf{c} \mapsto \overline{\mathbf{c}}$,
defines a real structure on the genus zero Deligne-Mumford moduli
space \dm{k}. \dm{k} is a real projective variety.
\end{lemma}
Proof. It is well-known that the tangent space
$T_{\mathbf{c}}\dm{k}$ at $\mathbf{c}$ is
$Ext^{1}(\Omega^{1}_{C}(a_{1}+ \ldots + a_{k}), \mathcal{O}_{C})$.
Thus,

\begin{gather}
T_{\mathbf{c}}\dm{k}  \notag \\
= \bigoplus_{i=1}^{l} H^{1}( \CP^{1}_{i}, T \CP^{1}_{i}(- \sum_{j}
q(a_{j,i}) - \sum_{\alpha,\beta} g_{\alpha}^{\beta,i})) \oplus
\bigoplus_{i=1}^{r} T_{g_{i}^{1}} \CP^{1}_{q(g_{i}^{1})}
\otimes T_{g_{i}^{2}} \CP^{1}_{q(g_{i}^{2})}, \label{e;cla} \\
\mbox{where $a_{j,i} \in \{ a_{1}, \ldots, a_{k} \}$ such that
$q(a_{j,i}) \in \CP^{1}_{i}$, and} \notag \\
\mbox{ $g_{\alpha}^{\beta,i} \in \pi^{-1}(g_{\alpha}) \bigcap
\CP^{1}_{i}$, where
 $g_{\alpha} \in \{g_{1},
\ldots, g_{r} \}$}     \notag\\
= \bigoplus_{i=1}^{k} T_{a_{i}}\CP^{1}_{q(a_{i})} \oplus
\bigoplus_{i=1,\ldots,r}^{j=1,2} T_{g_{i}^{j}}\CP^{1}_{q(g_{i}^{j})}
\ominus (\bigoplus_{i=1}^{l} H^{0}(\CP_{i}^{1}, T \CP^{1}_{i}))
\notag \\
 \oplus \bigoplus_{i=1}^{r} T_{g_{i}^{1}}
\CP^{1}_{q(g_{i}^{1})} \otimes T_{g_{i}^{2}} \CP^{1}_{q(g_{i}^{2})}.
\label{e;mordern}
\end{gather}

(\ref{e;cla}) comes from the following local to global spectral
sequence (cf. p99, \cite{ham}):

\begin{gather}
0 \hspace{1mm} \rightarrow \hspace{1mm} H^{1}(C,
\underline{Ext}^{0}_{C}(\Omega^{1}_{C}( a_{1} + \ldots + a_{k}),
\mathcal{O}_{C})) \hspace{1mm} \rightarrow  \label{e;logl} \\
\rightarrow \hspace{1mm} Ext^{1}(\Omega^{1}_{C}(a_{1}+ \ldots +
a_{k}), \mathcal{O}_{C}) \hspace{1mm} \rightarrow \hspace{1mm}
H^{0}(C, \underline{Ext}^{1}_{C}(\Omega^{1}_{C}(a_{1} + \ldots +
a_{k}),\mathcal{O}_{C})) \hspace{1mm} \rightarrow \hspace{1mm} 0
\notag
\end{gather}
See \cite{giv} for some further details. Terms $\bigoplus_{i=1}^{k}
T_{a_{i}}\CP^{1}_{q(a_{i})}$, $\bigoplus_{i=1,\ldots,r}^{j=1,2}
T_{g_{i}^{j}}\CP^{1}_{q(g_{i}^{j})}$, $\bigoplus_{i=1}^{l}
H^{0}(\CP_{i}^{1}, T \CP^{1}_{i})$ in (\ref{e;mordern}) came from
the long exact sequence of sheaf cohomology induced from the
following short exact sequence of sheaves:
\[ 0 \rightarrow T \CP^{1}_{i}(- \sum_{j}
q(a_{j,i}) - \sum_{\alpha, \beta} g_{\alpha}^{\beta,i}) \rightarrow
T \CP^{1}_{i} \rightarrow \bigoplus_{j} T_{q(a_{j,i})} \CP^{1}_{i}
\oplus \bigoplus_{\alpha, \beta} T_{g_{\alpha}^{\beta,i}}
\CP^{1}_{i} \rightarrow 0
\]

Signs in front of them are from the $K$-group calculation by using
Proposition 2.11. in \cite{aam}. Note that the rank of
$\displaystyle H^{0}(\CP^{1}_{i},T \CP^{1}_{i}(- \sum_{j} q(a_{j,i})
- \sum_{\alpha,\beta} g_{\alpha}^{\beta,i}) )$ is zero due to the
stability condition.

 The tangent space
$T_{\overline{\textbf{c}}}\dm{k}$ at $\overline{\textbf{c}}$ is:

\begin{gather}
 \bigoplus_{i=1}^{k} T_{s(a_{i})}\CP^{1}_{q(a_{i})} \oplus
\bigoplus_{i=1,\ldots,r}^{j=1,2}
T_{s(g_{i}^{j})}\CP^{1}_{q(g_{i}^{j})} \label{e;an} \\
 \ominus
(\bigoplus_{i=1}^{l} H^{0}(\CP_{i}^{1}, T \CP^{1}_{i})) \oplus
\bigoplus_{i=1}^{r} T_{s(g_{i}^{1})} \CP^{1}_{q(g_{i}^{1})} \otimes
T_{s(g_{i}^{2})} \CP^{1}_{q(g_{i}^{2})}. \notag
\end{gather}

The actual expression of the tangent space splitting depends on the
pointed curve representing the point $\mathbf{c}$. However, one can
observe the following. Let $(C',
\sigma(a_{1}),\ldots,\sigma(a_{k}))$ represent the point
$\mathbf{c}$ in $\dm{k}$, where  $\sigma$ is an element in
Aut($\CP^1$). Then, $(\overline{C'}, s \circ \sigma (a_{1}), \ldots,
s \circ \sigma(a_{k}))$ represents $\overline{\mathbf{c}}$.

\par Let $v$ be an
element in $ H^{0}(\CP_{i}^{1}, T \CP^{1}_{i})$. Let's denote $v
\mid_{x}$ be the value of $v$ at $x$. Then, $\bar{v}$ defined by
$\bar{v} \mid_{s(x)} := ds(v \mid_{x})$ is an element in
$H^{0}(\CP_{i}^{1}, T \CP^{1}_{i})$. The differential $ds$ of the
real structure $s$ induces the anti-holomorphic involution, $v
\mapsto \bar{v}$, on $H^{0}(\CP_{i}^{1}, T \CP^{1}_{i})$.
 The anti-holomorphic involutions on other components in
(\ref{e;mordern}) to (\ref{e;an}) are obviously induced by the
differential $ds$ on each component. Thus, the differential
$dI\mid_{\textbf{c}} : (\ref{e;mordern}) \mapsto (\ref{e;an})$ at
$\textbf{c}$ is an anti-holomorphic involution. $I$ is an
anti-holomorphic involution. \hfill $\Box$ \vspace{1mm}

\begin{remark}\label{r;demu}
\emph{We can also prove the Lemma\ref{L:D-M} as follows. The
Deligne-Mumford moduli space is originally defined over $\Z$. See
\cite[III.3]{mann}. So, it is defined over any field. The
$\C$-scheme Deligne-Mumford moduli space can be obtained by a scalar
extension from the $\R$-scheme Deligne-Mumford moduli space. That
is, the $\C$-scheme Deligne-Mumford moduli space is a
complexification $\dm{k}^{\R} \times_{\R}\C$ of the $\R$-scheme
Deligne-Mumford moduli space $\dm{k}^{\R}$. Thus, it has a canonical
anti-holomorphic involution. See \cite[p4, (1.4) Proposition]{sil}.
It is easily seen that the anti-holomorphic involution in Lemma
\ref{L:D-M} is identical to the corresponding canonical involution
in this Remark. }
\end{remark}

We calculate the tangent space on the moduli space $\md{X}$ of
$k$-pointed genus zero stable maps. Theorem \ref{p;tansp} is proven
in symplectic category by taking the different methods of
calculations in \cite{kwon3} when $f$ is an immersion on each
irreducible component. Intuitive interpretations of the
calculational results are seen in \cite{kwon3}.

\begin{theorem} \label{p;tansp}
Let $X$ be a convex projective variety. \\
 Let $\mathbf{f}:=[(f,C, a_{1}, \ldots, a_{k})]$ be a point in
$\md{X}$ such that $\beta$ is non-trivial. Let $\pi : \widetilde{C}
:= \CP^{1}_{1} \cup \ldots \cup \CP^{1}_{l} \rightarrow C$ be a
normalization map, where $\CP^{1}_{m}$ is biholomorphic to
$\CP^{1}$. Let $g_{1}, \ldots, g_{r}$ be singular points on $C$,
$r:= l-1$. Let's denote elements in $\pi^{-1}(g_{n})$ by $g_{n}^{1},
g_{n}^{2}$. Let $N_{m}$ be the normal sheaf induced from a morphism
$df_{m}:T\CP^{1}_{m}
\rightarrow T X$.\\
(i) Suppose $f_{m}:= f \mid_{\CP^{1}_{m}}$ is non-trivial, $m=1,
\ldots, l$. Then, the tangent space $T_{\mathbf{f}}\md{X}$ at
$\mathbf{f}$ is
\begin{gather}
\bigoplus_{m=1}^{l} H^{0}(\CP^{1}_{m}, N_{m}) \oplus \bigoplus_{i=1,
\ldots,k} T_{a_{i}}\CP^{1}_{q(a_{i})} \oplus
(\bigoplus_{n=1,\ldots,r} T_{g_{n}^{1}}\CP^{1}_{q(g_{n}^{1})}\otimes
T_{g_{n}^{2}} \CP^{1}_{q(g_{n}^{2})}) \oplus   \notag \\
 \oplus \bigoplus_{n=1, \ldots,r}^{j=1,2}
T_{g_{n}^{j}}\CP^{1}_{q(g_{n}^{j})} \ominus(\bigoplus_{n=1}^{r}
T_{f(g_{n})}X). \notag
\end{gather}

(ii) We may assume the following by reordering if necessary:\\
a. $f_{i}:= f \mid_{\CP^{1}_{i}}$ is trivial if $i=1, \ldots, m$,
and is non-trivial if $i=m+1, \ldots, l$.\\
b. If  $i=1, \ldots, h$, then, $g_{i}$ joins the irreducible
components on which the restriction of $f$ is non-trivial. If $i=
h+1, \ldots, r$, then $g_{i}$ joins the irreducible components such
that $f$ is trivial on one of the components or both components.
Then, the tangent space $T_{\mathbf{f}}\md{X}$ at $\mathbf{f}$ is
\begin{gather}
\bigoplus_{i=m+1}^{l} H^{0}(\CP^{1}_{i}, N_{i}) \oplus
\bigoplus_{i=1, \ldots,k} T_{a_{i}}\CP^{1}_{q(a_{i})} \oplus
(\bigoplus_{i=1,\ldots,r} T_{g_{i}^{1}}\CP^{1}_{q(g_{i}^{1})}\otimes
T_{g_{i}^{2}} \CP^{1}_{q(g_{i}^{2})}) \oplus   \notag \\
 \oplus \bigoplus_{i=1, \ldots,r}^{j=1,2}
T_{g_{i}^{j}}\CP^{1}_{q(g_{i}^{j})} \ominus(\bigoplus_{i=1}^{h}
T_{f(g_{i})}X) \ominus \bigoplus_{i=1}^{m} H^{0}(\CP^{1}_{i}, T
\CP^{1}_{i}). \notag
\end{gather}

\end{theorem}

\vspace{2mm}

Proof. We will use the following index notations throughout the
proof of (i):

\begin{tabular}{ll}
 $i= 1, \ldots, k$ & the index for marked points\\
$m \hspace{2mm} \text{or}\hspace{2mm} m'=1, \ldots, l$ & the index for irreducible components\\
 $n =
1, \ldots,  l-1$ & the index for gluing points\\
 $j= 1,2$ & the upper
index (with the lower index $n$) \\
& for pregluing points in $\pi^{-1}(g_{i})$.
\end{tabular}

\vspace{2mm}

To make the notations simpler, we will use the following notations:
\begin{center}
\begin{tabular}{ll}
$a_{i}(m)$ & $= a_{i}$ \hspace{2mm} if $a_{i} \in \CP^{1}_{m}$ \\ &
$=\emptyset$ \hspace{2mm} if $a_{i} \notin \CP^{1}_{m}$\\
 $g_{n}^{j}(m)$
 &$= g^{j}_{n} $\hspace{2mm} if \hspace{1mm}
$g^{j}_{n} \in \CP^{1}_{m}$ \\
& $ = \emptyset$ \hspace{2mm} if \hspace{1mm} $g_{n}^{j} \notin
\CP^{1}_{m}$ \\
$T_{\emptyset} \CP^{1}_{m} $ & $ := \emptyset$
\end{tabular}
\end{center}

For example, if $a_{1}(1) = \emptyset$, $a_{2}(1) = a_{2}$, then
$a_{1}(1) + a_{2}(1) = a_{2}$ and $T_{a_{1}(1)} \CP^{1}_{1}
\bigoplus T_{a_{2}(1)} \CP^{1}_{1} = T_{a_{2}} \CP^{1}_{1}$.
Obviously, $\sum_{i,m} a_{i}(m) = \sum_{i} a_{i}$.\\

As a convention, if we don't specify the range of the indices, eg.,
$a_{i}$, $i=1,2$, then we always consider all possible indexes. That
is, $a_{i}$ means
$a_{i}$, $i= 1, \ldots, k$.\\

The tangent space at $\textbf{f}$ is the hyperext group
$Ext^{1}(f^{*}\Omega^{1}_{X} \rightarrow \Omega^{1}_{C}(a_{1}+
\ldots + a_{k}), \mathcal{O}_{C})$.

 From \hfill the \hfill long
\hfill exact \hfill sequence \hfill associated \hfill with \hfill
the \hfill
hyperext \hfill group \\
$Ext^{1}(f^{*}\Omega^{1}_{X} \rightarrow \Omega^{1}_{C}(a_{1}+
\ldots + a_{k}), \mathcal{O}_{C})$: cf. \cite[ p285]{cak}

\begin{multline*}
0 \hspace{1mm} \rightarrow \hspace{1mm} Hom(\Omega_{C}^{1}(a_{1} +
\ldots + a_{k}), \mathcal{O}_{C}) \hspace{1mm} \rightarrow
\hspace{1mm}
H^{0}(C, f^{*}T X) \rightarrow\\
\rightarrow Ext^{1}(f^{*}\Omega^{1}_{X} \rightarrow
\Omega^{1}_{C}(a_{1} + \ldots + a_{k}), \mathcal{O}_{C})
\rightarrow\\
\rightarrow Ext^{1}(\Omega^{1}_{C}(a_{1} + \ldots +
a_{k}),\mathcal{O}_{C}) \rightarrow 0,
\end{multline*}

 we get the following tangent
space splitting at $\textbf{f}$:

\begin{equation}\label{e;ta}
\ominus Hom(\Omega_{C}^{1}(a_{1} + \ldots + a_{k}), \mathcal{O}_{C})
\oplus H^{0}(C, f^{*}T X) \oplus Ext^{1}(\Omega^{1}_{C}(a_{1} +
\ldots + a_{k}),\mathcal{O}_{C}).
\end{equation}

\vspace{2mm}

We will calculate each term's splitting first.

The standard fact we will use in the following calculations is
$Hom(\Omega_{C}^{1}(a_{1} + \ldots + a_{k}), \mathcal{O}_{C})$,
$\underline{Ext}^{0}(\Omega^{1}_{C}(a_{1}+ \ldots a_{k}),
\mathcal{O}_{C})$ is the sheaf of derivations gotten by the
pushforward of the sheaf of vector fields on $\widetilde{C}:=
\CP^{1}_{1} \cup \ldots \cup \CP^{1}_{l}$ vanishing at the inverse
images $g^{j}_{n}$ of the node in $C$
and the marked points $a_{i}$. cf. \cite[p100]{ham}\\

\par For $\ominus Hom(\Omega_{C}^{1}(a_{1} + \ldots + a_{k}), \mathcal{O}_{C})$ term,
we use the short exact sequences of sheaves:

\begin{multline}
0 \rightarrow T\CP^{1}_{m}(-\sum_{i} a_{i}(m) - \sum_{j,n}
g^{j}_{n}(m) ) \rightarrow T\CP^{1}_{m} \rightarrow \\
\rightarrow \bigoplus_{i}T_{a_{i}(m)}\CP^{1}_{m} \oplus
\bigoplus_{j,n} T_{g^{j}_{n}(m)}\CP^{1}_{m} \rightarrow 0
\label{e;hom1}
\end{multline}
\vspace{2mm}

to get the following K-group equation
\begin{align*}
  & Hom(\Omega_{C}^{1}(\sum_{i} a_{i} ), \mathcal{O}_{C})\\
= & H^{0}(C, TC(- \sum_{i} a_{i} ))\\
= & H^{0}(C, \pi_{*}( T\widetilde{C}(-\sum_{i} a_{i} - \sum_{j,n}
g^{j}_{n})))
\end{align*}
\begin{align*} = & H^{0}(\widetilde{C},  T \widetilde{C}(- \sum_{i}
a_{i} -
\sum_{j,n}g^{j}_{n}))\\
 = & \bigoplus_{m} H^{0}(\CP^{1}_{m},
T\CP^{1}_{m}(-\sum_{i} a_{i}(m)
- \sum_{j,n} g^{j}_{n}(m)))\\
= & \bigoplus_{m} \textbf{[}H^{0}(\CP^{1}_{m}, T \CP^{1}_{m})
 \ominus (\bigoplus_{i} T_{a_{i}(m)}\CP^{1}_{m})\\
& \ominus (\bigoplus_{j,n} T_{g^{j}_{n}(m)}\CP^{1}_{m}) \oplus
H^{1}(\CP^{1}_{m}, T \CP^{1}_{m}(-\sum_{i} a_{i}(m) - \sum_{j,n}
g^{j}_{n}(m)))\textbf{]}
\end{align*}

\noindent by (~\ref{e;hom1}).\\

\par For $H^{0}(C, f^{*}T X)$, we use the short exact sequence of
sheaves

\begin{center} $0 \rightarrow f^* T X
\rightarrow \bigoplus_{m} f_{m}^{*} TX \rightarrow \bigoplus_{n}
T_{f(g_{n})} X \rightarrow 0$,
\end{center}
\vspace{2mm}

\noindent to get the K-group equation

\begin{center}
$H^{0}(C, f^{*}T X) = \bigoplus_{m} H^{0}(\CP^{1}_{m}, f_{m}^{*} T
X) \ominus \bigoplus_{n} T_{f(g_{n})}X$,
\end{center}
\vspace{2mm}

\noindent because $H^{1}(C, f^{*}T X)$ vanishes by Lemma 10 in
~\cite{fap}.\\

\par For $Ext^{1}(\Omega^{1}_{C}(a_{1} + \ldots + a_{k}),\mathcal{O}_{C})$,
we use the exact sequence from the local to global spectral sequence
in (\ref{e;logl}) to get

\begin{align*}
& Ext^{1}(\Omega^{1}_{C}(a_{1} + \ldots + a_{k}),\mathcal{O}_{C})\\
= &  H^{1}(C, \underline{Ext}^{0}(\Omega_{C}(\sum_{i} a_{i}),
\mathcal{O}_{C})) \oplus H^{0}(C,
\underline{Ext}^{1}(\Omega_{C}(\sum_{i} a_{i} ), \mathcal{O}_{C})
\end{align*}
\begin{align*}
 =
& H^{1}(C, \pi_{*}( T \widetilde{C}(-\sum_{i} a_{i} - \sum_{i,n}
g^{j}_{n} )))
 \oplus H^{0}(C, \underline{Ext}^{1}(\Omega_{C}(\sum_{i}a_{i}
), \mathcal{O}_{C}))\\
= & \bigoplus_{m} \textbf{[} H^{1}(\CP^{1}_{m},
T\CP^{1}_{m}(-\sum_{i} a_{i}(m) -
\sum_{n,j}g^{j}_{n}(m)))\textbf{]}\\
& \oplus
 \bigoplus_{m,m',n}[(T_{g_{n}^{1}(m)} \CP^{1}_{m})
  \otimes (T_{g_{n}^{2}(m')} \CP^{1}_{m'})].
\end{align*}

\par The result follows by putting all terms to (\ref{e;ta}) and simplify further
by $K$-group calculations with the long exact sequence of  sheaf
cohomology induced from the following the short exact sequence of sheaves:\\

\begin{center}
$0 \rightarrow T \CP^{1}_{m} \rightarrow f^{*}_{m} T X \rightarrow
N_{m} \rightarrow 0$
\end{center}

\par The proof of (ii) is very similar to the proof of (i). Use
that if $f_{\alpha}$ is trivial, then $H^{0}(\CP^{1}_{\alpha},
f_{\alpha}^{*} TX)$ is isomorphic to $T_{f_{\alpha}
(\CP^{1}_{\alpha})}X$.  \hfill $\Box$

\vspace{2mm}

\par For the proof of Theorem \ref{t;maint}, it is enough to show that the involution defined by
$[(f, C, a_{1}, \ldots, a_{k})]  \mapsto [(\overline{f},
\overline{C}, s(a_{1}), \ldots, s(a_{k}))]$ on the moduli space
$\md{X}$ is an anti-holomorphic involution( cf.\cite[p4, (1.4)
Proposition]{sil}). $\md{X}$ is a normal projective variety. It has
orbifold singularities. So, we show that the defined involution is
an anti-holomorphic involution with local chart before the local
quotient by a finite group action and then show there is a canonical
conjugate group action on the conjugate local chart around the
conjugate point.

\begin{theorem}\label{t;maint}
Let $X$ be a convex real projective variety. Then, the moduli space
$\md{X}$ of stable maps is a real projective variety whose real
structure comes from the involution defined by $[(f, C, a_{1},
\ldots, a_{k})] \mapsto [( t \circ f \circ s, \overline{C},
s(a_{1}), \ldots, s(a_{k}))]$, where the notations $C$ and
$\overline{C}$ are the same as in Lemma \ref{L:D-M}.
\end{theorem}

Proof. Let $\overline{\mathbf{f}}$ be a point in $\md{X}$
represented by $(\overline{f}, \overline{C}, s(a_{1}), \ldots,
s(a_{k}))$, where $\overline{f}(z) := t \circ f\circ s(z)$. Let $H:
\md{X} \rightarrow  \md{X}$ be the involution defined by $\mathbf{f}
\mapsto \overline{\mathbf{f}}$. The Theorem follows if we show that
$H$ is an anti-holomorphic involution.

\par Let \hfill $\beta$ \hfill is \hfill non-trivial. \hfill
Let's \hfill suppose \hfill that \hfill $f$ \hfill is \hfill non-trivial \hfill on
\hfill every \hfill component. \hfill From \hfill
Theorem\ref{p;tansp} (i), \hfill we \hfill know \hfill that \hfill
the  \hfill tangent \hfill space \hfill
at\\
$\mathbf{f}:=[(f,C,a_{1}, \ldots, a_{k})]$ is

\begin{gather} \bigoplus_{i=1}^{l} H^{0}(\CP^{1}_{i}, N_{i})
\oplus \bigoplus_{i=1, \ldots,k} T_{a_{i}}\CP^{1}_{q(a_{i})} \oplus
(\bigoplus_{i=1,\ldots,r} T_{g_{i}^{1}}\CP^{1}_{q(g_{i}^{1})}\otimes
T_{g_{i}^{2}} \CP^{1}_{q(g_{i}^{2})}) \oplus   \notag \\
 \oplus \bigoplus_{i=1, \ldots,r}^{j=1,2}
T_{g_{i}^{j}}\CP^{1}_{q(g_{i}^{j})} \ominus(\bigoplus_{i=1}^{r}
T_{f(g_{i})}X), \notag
\end{gather}

\par Let $\overline{N}_{i}$ be the normal sheaf induced from the
morphism $d\overline{f}_{i}: T \CP^{1}_{i} \rightarrow TX$, where
$\overline{f}_{i}(z):= t \circ f_{i} \circ s(z)$. The tangent space
at $\overline{\textbf{f}}$ is:

\begin{gather} \bigoplus_{i=1}^{l} H^{0}(\CP^{1}_{i}, \overline{N}_{i})
\oplus \bigoplus_{i=1, \ldots,k} T_{s(a_{i})}\CP^{1}_{q(a_{i})}
\oplus (\bigoplus_{i=1,\ldots,r}
T_{s(g_{i}^{1})}\CP^{1}_{q(g_{i}^{1})}\otimes T_{s(g_{i}^{2})}
\CP^{1}_{q(g_{i}^{2})}) \oplus   \notag \end{gather}
\begin{gather}
 \oplus \bigoplus_{i=1, \ldots,r}^{j=1,2}
T_{s(g_{i}^{j})}\CP^{1}_{q(g_{i}^{j})} \ominus(\bigoplus_{i=1}^{r}
T_{t \circ f(g_{i})}X). \notag
\end{gather}

Each term in the tangent space splitting at $\mathbf{f}$,
$\overline{\mathbf{f}}$ is a complex vector space.

\begin{gather}
\bigoplus_{i=1, \ldots,k} T_{a_{i}}\CP^{1}_{q(a_{i})}
\stackrel{dH}{\mapsto} \bigoplus_{i=1, \ldots,l}
T_{s(a_{i})}\CP^{1}_{q(a_{i})} \label{e;fi} \\
\bigoplus_{i=1,\ldots,r} T_{g_{i}^{1}}\CP^{1}_{q(g_{i}^{1})}\otimes
T_{g_{i}^{2}} \CP^{1}_{q(g_{i}^{2})} \stackrel{dH}{\mapsto}
\bigoplus_{i=1,\ldots,r}
T_{s(g_{i}^{1})}\CP^{1}_{q(g_{i}^{1})}\otimes
T_{s(g_{i}^{2})} \CP^{1}_{q(g_{i}^{2})} \label{e;se}\\
\bigoplus_{i=1, \ldots,r}^{j=1,2}
T_{g_{i}^{j}}\CP^{1}_{q(g_{i}^{j})} \stackrel{dH}{\mapsto}
\bigoplus_{i=1, \ldots,r}^{j=1,2}
T_{s(g_{i}^{j})}\CP^{1}_{q(g_{i}^{j})} \label{e;th}
\end{gather}

It is obvious that $dH$ in (\ref{e;fi}), (\ref{e;se}), (\ref{e;th})
is the anti-holomorphic involution induced by the real structure of
a complex conjugation map on $\CP^1$.
\begin{gather}
\bigoplus_{i=1}^{r} T_{f(g_{i})}X \stackrel{dH}{\mapsto}
\bigoplus_{i=1}^{r} T_{t \circ f(g_{i})}X \label{e;targ}
\end{gather}

Clearly, $dH$ in (\ref{e;targ}) is the anti-holomorphic involution
induced by the real structure $t$ on the target space $X$.

\begin{gather} \bigoplus_{i=1}^{l} H^{0}(\CP^{1}_{i}, N_{i})
\stackrel{dH}{\mapsto}  \bigoplus_{i=1}^{l} H^{0}(\CP^{1}_{i},
\overline{N}_{i})    \label{e;fin}
\end{gather}

Each normal sheaf $N_{i}$ is the direct sum of the locally free
sheaf $NB_{i}$ and skyscraper sheaves supported by critical points
of $f$. Similar to the case (\ref{e;targ}), the restriction of $dH$
to each skyscraper sheaf is the anti-holomorphic involution induced
by the real structure $t$ on the target space $X$. Normal bundles
$NB_{i}$, $\overline{NB}_{i}$ split into line bundles on $\CP^1$ by
splitting principle. By considering Weil divisors characterizing
each line bundle and the definition of $\overline{f}_{i}$, one can
check that $\overline{NB}_{i}$ is a conjugate bundle for the bundle
$NB_{i}$. Thus, the restriction of $dH$ to $H^{0}(\CP^{1}_{i},
NB_{i})$ is an anti-holomorphic involution.

\par The general case considered in Theorem\ref{p;tansp} (ii) can be
easily proven by repeating the same arguments we did above, that is,
by checking the componentwise anti-holomorphicity of $dH$.

\par Let's assume that $(f', C', b_{1}, \ldots, b_{k})$ represents
$\mathbf{f}$. Then, there exists the element $\sigma \in$
Aut($\CP^1$) such that $f= f' \circ \sigma$ and $b_{i} =
\sigma(a_{i})$. $\mathbf{\overline{f}}$ is represented by
$(\overline{f},\overline{C}, s(a_{1}), \ldots, s(a_{k}))$. Note that
$\overline{f}' \circ \overline{\sigma} = t \circ f' \circ s \circ s
\circ \sigma \circ s = t \circ f' \circ \sigma \circ s = t \circ f
\circ s = \overline{f}$, $s(b_{i}) = s \circ \sigma \circ
s(s(a_{i}))$. One can check that $s \circ \sigma \circ s$ is also an
element in Aut($\CP^1$). Thus, it shows that $( t \circ f' \circ s,
\overline{C}', s(b_{1}), \ldots,s(b_{k}))$ represents
$\overline{\textbf{f}}$. This implies that the map $H$ is an
anti-holomorphic involution on the local charts, independent of the
actual choice of the chosen pointed stable map representing
$\textbf{f}$.
\par Let $G$ be a finite group acting on the local chart $\mathcal{O}_{\textbf{h}}$
around $\textbf{h}$ such that $ G \times \mathcal{O}_{\textbf{h}}
\rightarrow \mathcal{O}_{\textbf{h}}, (g, \textbf{f}) \mapsto g
\cdot \textbf{f}$. Let
$\overline{\mathcal{O}}_{\overline{\textbf{h}}}$ be the conjugate
local chart around $\overline{\textbf{h}}$, i.e.,
$\overline{\mathcal{O}}_{\overline{\textbf{h}}} :=
\{\overline{\textbf{f}} \mid \textbf{f} \in \mathcal{O}_{\textbf{h}}
\}$. Then, there is a canonical conjugate group action of $G$
defined by $g \cdot \overline{\textbf{f}}:= \overline{g \cdot
\textbf{f}}$
\par Thus, $H$ is an anti-holomorphic involution on the
orbifold $\md{X}$.
\par Let's assume that $\beta$ is trivial. Then, the moduli space
$\md{X}$ is isomorphic to the complex manifold $\dm{k} \times X$. It
is obvious the defined involution is an anti-holomorphic involution.
\hfill $\Box$

\begin{remark}\label{R;kollar}
\emph{Araujo - Kollar constructed the moduli space of stable maps on
any Noetherian scheme in \cite[sec.10]{aak}. It is interesting to
see whether the variety gotten by the complexification of the moduli
space defined over $\R$ is isomorphic to the moduli space of stable
maps defined over $\C$ or not.  Note that the complexification as a
variety doesn't need to have any meaning as a moduli space.
Nevertheless, an element in the real part of the variety gotten by
the complexification uniquely corresponds to a real point in the
moduli space constructed over $\R$.}
\par \emph{ Let's consider the degree 2 maps from $\CP^1$ to $\CP^1(:= \C \cup \{ \infty \})$,
defined by $z \mapsto z^{2}$ and $z \mapsto - z^{2}$. Then, they are
represented by two distinct real points in the moduli space
constructed over $\R$. But they are represented by one point in the
real part of the complex moduli space $\overline{M}_{0}(\CP^1, 2
\cdot \text{[line]})$ because they are equivalent maps. Thus, the
complexification of the moduli space of degree 2 maps over $\R$
cannot be isomorphic to the complex moduli space
$\overline{M}_{0}(\CP^1, 2 \cdot \text{[line]})$.}
\par \emph{Due to
the differences of the categories and equivalence relations, the
complexification of the moduli space defined over $\R$ as a variety
is not always isomorphic to the moduli space defined over $\C$.
Recall Remark \ref{r;demu} for the Deligne-Mumford moduli space. }
\end{remark}

\section{Real Properties of the Moduli space}\label{s;property}

\begin{proposition}\label{p;ev}
The $i$-th evaluation map $ev_{i}$ is a real map.
\end{proposition}

Proof. It is enough to show that $ev_{i}$ commutes with the real
structure $H$ on $\md{X}$, $t$ on $X$. See \cite[p107,
4.7.(c)]{hart}. Let $\mathbf{f}:= [(f,C, a_{1}, \ldots, a_{k})]$.
Then, $H(\mathbf{f}) :=[(t \circ f \circ s, s(a_{1}), \ldots,
s(a_{k}))]$. Thus, $t\circ f \circ s(s(a_{i}))
=ev_{i}(H(\mathbf{f})) = t(ev_{i}(\mathbf{f})) = t(f(a_{i}))$. This
commutation relation is independent of the pointed stable map
representing the point $\mathbf{f} \in \md{X}$.
Thus, the Proposition follows. \hfill $\Box$ \\
\par One can show the following Proposition as we proved Proposition
\ref{p;ev}. Its proof is left to readers.

\begin{proposition}\label{p;for}
The forgetful maps, $\md{X} \rightarrow \overline{M}_{k-1}(X,
\beta)$, $\md{X} \rightarrow \dm{k}$, are real maps.
\end{proposition}

\begin{proposition}\label{p;realp} Let $\CP^n$ have the real structure from a
complex conjugation involution. Let $X$ be a real projective variety
such that the imbedding $i$ which decides the real structure of $X$
intersects with the real part $\RP^n$ of $\CP^n$. If $k \geq 3$,
then the real part of $M_{k}(X, \beta)$( the locus in $\md{X}$ whose
domain curve is smooth ) consists of stable maps $[(f, \CP^1, a_{1},
\ldots , a_{k})]$ such that $a_{i}$ is in $\RP^1 (\subset \CP^1)$
and $f$ is a real map, i.e., $f \circ s = t \circ f$.
\end{proposition}

Proof. It is well-known that the real part of the Deligne-Mumford
moduli space $M_{k}$ (before the compactification) consists of
curves whose marked points are on the real part $\RP^{1}$ of the
domain curve $\CP^{1}$. See \cite[sec.2.3]{gam}. If $\mathbf{f} :=
[(f, \CP^{1}, a_{1}, \ldots, a_{k})]$ represents a point in the real
part of $M_{k}(X, \beta)$, then Proposition \ref{p;for} asserts that
$\mathbf{f}$ is represented by a map $(f, \CP^1, a_{1}, \ldots,
a_{k})$, where $a_{i}$ is on the real part $\RP^1$ in $\CP^1$. And
there exists $\sigma \in$ Aut($\CP^1$) such that $f = t \circ f
\circ s \circ \sigma$, $a_{i} = \sigma(a_{i})$. Since $k \geq 3$,
$\sigma$ is an identity map. It
shows that $f$ is a real map. \hfill $\Box$\\

\par Let's assume that the target space $X$ has the same real
structure stated in Proposition \ref{p;realp}. Then, of course, the
type of stable maps described in Proposition \ref{p;realp} are in
the real part of $M_{k}(X, \beta)$. But, the real part analysis of
the whole moduli space $\md{X}$ for all $k$ is subtle. Let's
illustrate with some examples.
\begin{itemize}
\item  If $k=0$, then some non-real maps are in the real part of the
moduli space. For example, the non-real map $f: \CP^1 \rightarrow
\CP^2$ whose image curve is represented by the equation $x^{2} +
y^{2} + z^{2}$ in $\CP^2$ is in the real part of the moduli space.
\item On $\overline{M}_{k}(X, \beta) \setminus M_{k}(X, \beta)$, the
stable maps all of whose gluing points are in the real part of the
domain curve and all of whose marked points are in the real part of
each irreducible component are in the real part of the moduli space.
\item The gluing points in the reducible domain curve don't have
to be in the real part of the domain curve. Let $[(f, C, a_{1})]$ be
the element in $\overline{M}_{1}(\CP^1, 2 \cdot \text{[line]})$ such
that \par $\ast$ the normalization $\widetilde{C}$ of the domain
curve $C$ is $\CP^{1}_{0} \cup \CP^{1}_{1} \cup \CP^{1}_{2}$, where
$\CP^{1}_{i} \cong \CP^{1} \cong \C \cup \{ \infty \}$
\par $\ast$ the point $0 \in \CP^{1}_{1}$ is glued to the point $i
\in \CP^{1}_{0}$ and the point $ 0 \in \CP^{1}_{2}$ is glued to the
point $-i \in \CP^{1}_{0}$ \par $\ast$ $f \mid_{\CP^{1}_{0}} = 0$,
$f \mid_{\CP^{1}_{1}} = f \mid_{\CP^{1}_{2}} =$ identity map
\end{itemize}

The last example was given by Pierre Deligne to the author.

\vspace{1mm}

\noindent \textbf{Acknowledgments.} I thank my advisor Selman
Akbulut for suggesting the subject and supports when I was a
graduate student. I also thank the referee for many helpful
suggestions.

\vspace{1cm}

\begin{center}
Oklahoma State University \\
Department of Mathematics \\
401 Mathematical Sciences, Stillwater, OK 74078-1058\\
U.S.A.\\
 seonkwon@math.okstate.edu
\end{center}

\end{document}